\documentclass[11pt,reqno]{amsart}
\usepackage{graphicx}

\numberwithin{equation}{section}
\newtheorem{theorem}{Theorem}[section]
\newtheorem{lemma}[theorem]{Lemma}

\newtheorem{condition}[theorem]{Condition}

\theoremstyle{definition}

\newtheorem{remark}[theorem]{Remark}

\newenvironment{romenumerate}{\begin{enumerate}
 }{\end{enumerate}}

\newcounter{thmenumerate}

\newcommand{\refT}[1]{Theorem~\ref{#1}}

\newcommand{\refL}[1]{Lemma~\ref{#1}}
\newcommand{\refCN}[1]{Condition~\ref{#1}}

\newcommand{\refS}[1]{Section~\ref{#1}}
\newcommand{\refSS}[1]{Subsection~\ref{#1}}

\newcommand{\refand}[2]{\ref{#1} and~\ref{#2}}

\newcommand{\qedkorr}{\vskip-\baselineskip} 

\begingroup
  \divide\count255 by 60
  \multiply\count255 by -60
  \advance\count255 by \time
  \ifnum \count255 < 10 \xdef\klockan{\the\count1.0\the\count255}
  \else\xdef\klockan{\the\count1.\the\count255}\fi
\endgroup

\newcommand\webcite[1]{\hfil\penalty0\texttt{\def~{\~{}}#1}\hfill\hfill}

\newcommand\set[1]{\ensuremath{\{#1\}}}

\newcommand\bigpar[1]{\bigl(#1\bigr)}
\newcommand\Bigpar[1]{\Bigl(#1\Bigr)}

\newcommand\lrpar[1]{\left(#1\right)}
\def\rompar(#1){\textup(#1\textup)}    

\newcommand\ntoo{\ensuremath{{n\to\infty}}}

\newcommand\ktoo{\ensuremath{{k\to\infty}}}

\newcommand\iid{i.i.d.\spacefactor=1000}     
\newcommand\ie{i.e.\spacefactor=1000}
\newcommand\eg{e.g.\spacefactor=1000}
\newcommand\viz{viz.\spacefactor=1000}
\newcommand\cf{cf.\spacefactor=1000}
\newcommand{\as}{a.s.\spacefactor=1000}

\newcommand\whp{{whp}}

\newcommand\dto{\overset{\mathrm{d}}{\to}}
\newcommand\pto{\overset{\mathrm{p}}{\to}}

\renewcommand\={:=}
\newcommand\bbR{\mathbb R}

\newcommand\bbN{\mathbb N}

\newcommand\E{\operatorname{\mathbb E{}}}
\renewcommand\P{\operatorname{\mathbb P{}}}
\newcommand\Var{\operatorname{Var}}

\newcommand\Po{\operatorname{Po}}
\newcommand\Bi{\operatorname{Bi}}
\newcommand\MixBi{\operatorname{MixBi}}
\newcommand\Be{\operatorname{Be}}

\newcommand\sumin{\sum_{i=1}^n}

\newcommand\sumdoo{\sum_{d=1}^\infty}
\newcommand\ga{\alpha}

\newcommand\gd{\delta}

\newcommand\kk{\kappa}

\newcommand\go{\omega}

\newcommand\eps{\varepsilon}

\newcommand\cC{\mathcal C}

\newcommand\cL{{\mathcal L}}
\newcommand\tD{\tilde D}
\newcommand\tX{{\tilde X}}
\newcommand\tp{\tilde p}
\newcommand\tpi{\tilde \pi}

\newcommand\qi{^{-1}}

\newcommand\pij{p_{ij}}

\newcommand\gnx[1]{\ensuremath{G(n,#1)}}

\newcommand\gnk{\gnx{\kk}}

\newcommand\gnp{\ensuremath{G(n,p)}}
\newcommand\gnm{\ensuremath{G(n,m)}}

\newcommand\gnd{\ensuremath{G(n,(d_i)_1^n)}}
\newcommand\gndx{\ensuremath{G^*(n,(d_i)_1^n)}}
\newcommand\gndxa{\ensuremath{\gndx\acp}}
\newcommand\gnw{\ensuremath{G_W(n)}}
\newcommand\nn{^{(n)}}

\newcommand\bp{\ensuremath{\mathfrak{X}}}
\newcommand\bpd{\ensuremath{\mathfrak{X}^d}}
\newcommand\bpnu{\ensuremath{\mathfrak{X}_{\nu}}}
\newcommand\bpnud{\ensuremath{\mathfrak{X}_{\nu}^d}}

\newcommand\bpxd{\ensuremath{\mathfrak X^{(d)}}}
\newcommand\bpxdnu{\ensuremath{\mathfrak X^{(d)}_\nu}}

\newcommand\nx[1]{\ensuremath{N_{#1}}}
\newcommand\nk{\nx{k}}
\newcommand\ngek{\nx{\ge k}}
\newcommand\ngex[1]{\nx{\ge #1}}
\newcommand\ngegon{\ngex{\go(n)}}
\newcommand\ndk{\nx{k,d}}
\newcommand\ndgek{\ndgex{k}}
\newcommand\ndgex[1]{\nx{\ge #1,d}}
\newcommand\ndgegon{\ndgex{\go(n)}}

\newcommand\CS{Cauchy--Schwarz}

\newcommand{\vvq}{\mathsf{v}}
\newcommand{\vv}{_{\vvq}}
\newcommand{\vvp}{_{\vvq; p}}
\newcommand{\ux}{\mathsf{U}}
\newcommand{\uq}{^{\ux}}
\newcommand{\uv}{^{\ux}_{v}}
\newcommand{\uvp}{^{\ux}_{v; p}}
\newcommand{\ax}{\mathsf{A}}
\newcommand{\aq}{^{\ax}}
\newcommand{\ac}{^{\ax}_{c}}
\newcommand{\acp}{^{\ax}_{c; p}}

\newcommand{\vc}{(c)}
\newcommand{\eix}{\mathsf{E1}}
\newcommand{\ei}{^{\eix}}
\newcommand{\eiq}{^{\ei}_{\alpha}}
\newcommand{\eiqp}{^{\ei}_{\alpha; p}}
\newcommand{\eiix}{\mathsf{E2}}
\newcommand{\eii}{^{\eiix}}
\newcommand{\eiiq}{^{\eii}_{\alpha}}
\newcommand{\eiiqp}{^{\eii}_{\alpha ; p}}
\newcommand{\ejx}{\mathsf{E}j}
\newcommand{\ej}{^{\ejx}}
\newcommand{\ejq}{^{\ej}_{\alpha }}
\newcommand{\ejqp}{^{\ej}_{\alpha; p}}

\newcommand\xp{\bar p}
\newcommand{\nigd}{n^{1-\gd}}
\newcommand{\ngd}{n^{-\gd}}

\newcommand{\pidu}{\pi^{(d)}}
\newcommand{\pida}{\pi^{(d)}}

\begin{document}

\title[Graphs, epidemics and vaccination strategies]
{Graphs with specified degree distributions, simple epidemics
  and local vaccination strategies}
\author{Tom Britton, Svante Janson and Anders Martin-L\"of}

\date{January 25, 2007} 

\address{Tom Britton, Department of Mathematics, Stockholm University,
SE-10691 Stockholm, Sweden.} 
\email{tom.britton@math.su.se}

\address{Svante Janson, Department of
Mathematics, Uppsala University, P.O. Box 480, SE-75106 Uppsala, Sweden.}
\email{svante.janson@math.uu.se}

\address{Anders
Martin-L\"of, Department of Mathematics, Stockholm University,
SE-10691 Stockholm, Sweden.}
\email{andersml@math.su.se}



\begin{abstract}
Consider a random graph, having a pre-specified degree
distribution $F$ but other than that being uniformly distributed,
describing the social structure (friendship) in a large community. Suppose
one individual in the community is externally infected
by an infectious disease and that the disease has its course
by assuming that infected individuals infect their not yet infected friends
independently with probability $p$. For this situation the paper
determines $R_0$ and $\tau_0$, the basic reproduction number and
the asymptotic final size in case of a major outbreak. Further,
the paper looks at some different local vaccination strategies
where individuals are chosen randomly and vaccinated, or friends
of the selected individuals are vaccinated, prior to the
introduction of the disease. For the studied vaccination
strategies the paper determines $R_v$: the reproduction number, and
$\tau_v$: the asymptotic final proportion infected in case of a
major outbreak, after vaccinating a fraction $v$.
\end{abstract}

\maketitle

\section{Introduction}

Simple undirected random graphs can be used to describe the social
network in a large community (e.g.\ \cite{s00}), vertices corresponding to
individuals and edges to some type of social relation, from now on
denoted friendship. Given such a graph, a model for the spread of the
disease may be defined, where individuals at first are susceptible but may then
become infected by a friend. An infected individual has the
potential to spread the disease to its not yet infected friends before
it recovers and becomes immune. The final outbreak, both its size and who
gets infected, depends on properties of the social graph as well
as on properties of disease transmission. In order to prevent an outbreak it is
possible to vaccinate, or immunize in some other way, individuals
prior to arrival of the disease. Who and how many that are to be
vaccinated specifies the vaccination strategy.

The present paper studies questions arising from such modeling. In
particular, we consider random graphs where the degree distribution
(the number of friends) follows some pre-specified distribution
$F$, typically having heavy tails, but where the random graph $G$ is
otherwise uniformly distributed. The epidemic model is the
simplest possible model for a susceptible-infectious-removed (SIR)
disease (e.g.\ \cite{ab00}). One randomly selected individual is
initially externally infected. Any individual who
becomes infected infects each of his/her not yet infected friends
independently with probability $p$, and after that the individual
recovers and becomes immune, a state called removed. For this graph
and epidemic model we study
different vaccination strategies: the uniform strategy and the
acquaintance strategy \cite{chb03}. In both strategies
individuals are chosen randomly from the community. In the uniform
strategy the selected individuals are vaccinated and in the
acquaintance strategy a randomly chosen friend of the selected
individual is vaccinated. Both vaccination strategies are local in
the sense that the global social network need not be known in
order to perform the strategy. We also study a vaccination strategy
where, instead of selecting individuals at random, friendships are
selected and one or two of the corresponding friends get vaccinated.

As the population size $n$ tends to infinity, we prove that the
initial phase of the epidemic may be approximated by a suitable
branching process. The largest eigenvalue of the branching
process, often denoted $R_0$ and called the basic reproduction
number when applied to epidemics \cite{ab00},
determines whether a major outbreak can occur or not: if $R_0\le
1$ only minor outbreaks can occur whereas if
$R_0>1$ outbreaks of order $O(n)$ can also occur with positive
probability. In case of a
major outbreak the total number of individuals infected during the
outbreak, the final size, is shown to satisfy a law of large number.
The corresponding (random) proportion is shown to converge in
probability to a deterministic
limit $\tau_0$. Similar results are obtained when a vaccination
strategy with vaccination coverage $v$ has been performed prior to
disease introduction. In this situation the strategy-specific
reproduction number $R_v$, and the major outbreak size $\tau_v$, are
determined. From this it is possible to determine the
(strategy-specific) critical vaccination coverage $v_c$ which
determines the necessary 
proportion to vaccinate in order to surely prevent a major outbreak, so
$v_c=\inf_v \{v; R_v\le 1\}$.

Stochastic epidemic models on networks with pre-specified degree
distributios have
mainly been studied in the physics literature (e.g. \cite{mn00},
\cite{nsw01}, \cite{chb03}), Andersson \cite{a99} being one
exception. Some of the problems
studied in the present paper have been analysed before whereas
others have not, in particular the final size proportion $\tau_v$
as a function of $v$. Beside contributing with some new results
another aim of the paper is to give formal proofs to results which have
previously only been obtained heuristically.

The rest of the paper is structured as follows. In \refS{Smodels} we
define the models for the random graph, the epidemic and the
vaccination strategies. In \refS{Sresults} we present the main results,
motivate them with some heuristics and give some examples and
illustrations. The
proofs are given in Sections \refand{Sgw}{Spf}.

\section{Models}\label{Smodels}

\subsection{Graphs}\label{SSgraphs}

Let $G$ denote a random \emph{multigraph}, allowing for multiple edges and
loops, and let $n=|G|$ denote the number of
vertices of $G$, \ie{} the population size.
Later we shall consider limits as \ntoo. We define our random
multigraph as follows. Let $n \in \bbN$ and let
$(d_i)_1^n=(d_i\nn)_1^n$ be a sequence of
non-negative integers such that  $\sumin d_i$ is even.
We define a \emph{random multigraph with given degree
sequence $(d_i)_1^n$}, denoted by \gndx,
by the configuration model (see \eg{} \cite{bollobas}):
take a set of $d_i$ half-edges for each vertex
$i$, and combine the half-edges into pairs by a uniformly random
matching of the set of all half-edges.

Note that \gndx{} does not have exactly the uniform
distribution over all multigraphs with the given degree sequence;
there is a weight with a factor $1/j!$ for every edge of multiplicity $j$, and
a factor $1/2$ for every loop, see \cite[\S1]{JKLP}. However,
conditioned on the multigraph being a (simple) graph, we obtain a
uniformly distributed random graph with the given degree sequence,
which we denote by \gnd. It is also worth mentioning that
the distribution of \gndx{} is the same as the one obtained by sampling the
edges as ordered pairs of vertices uniformly with replacement, and
then conditioning on the vertex degrees being correct.

Let us write $2m\=\sumin d_i$, so that $m=m(n)$ is the number of
edges in the multigraph \gndx.
We assume that we are given
$(d_i)_1^n$ satisfying the following regularity conditions,
\cf{} Molloy and Reed \cite{MR1,MR2}.
\begin{condition}\label{C1}
For each $n$,  $(d_i)_1^n=(d_i\nn)_1^n$ is a sequence of non-negative
integers such that $\sumin d_i$ is even and, for some probability distribution
$(p_j)_{j=0}^\infty$ independent of $n$,
and with
$n_j\=\#\set{i:d_i=j}$,
\begin{romenumerate}
  \item
$n_j/n\to p_j$ for every $j\ge0$ as \ntoo;
\item
$\mu\=\sum_j j p_j\in(0,\infty)$;
\item
$2m/n\to\mu$ as \ntoo.
\item $p_2<1$.
\end{romenumerate}
\end{condition}

\begin{remark}\label{Runiform}
Note that $2m=\sum_i d_i = \sum_j j n_j$. Thus, \refCN{C1} implies
that the sum $\sum_j jn_j/n$ converges uniformly for $n\ge1$,
\ie
\begin{equation}\label{unif1}
  \lim_{J\to\infty} \sup_n \sum_{j>J} j n_j/n =0.
\end{equation}
Conversely, \eqref{unif1}
together with (i) and (ii) implies (iii).
(This follows from, \eg, \cite[Theorem 5.5.4]{Gut}, taking $X_n$ to be
the degree of a random vertex.)

Note that our condition is slightly weaker than the one in
Molloy and Reed \cite{MR1,MR2}; 
they also assume (in an equivalent formulation)
that if $\sum_j j^2 p_j<\infty$, then the sums $\sum_j j^2 n_j/n$
converge uniformly; moreover they assume that $j^2 n_j/n\to j^2 p_j$
uniformly.
\end{remark}

\refCN{C1} is all we need to study the random multigraph \gndx. In 
order to treat the random simple graph \gnd, which is our main model,
we need an additional assumption.
\begin{condition}\label{C2}
$\sum_i d_i^2=O(n)$.
\end{condition}
Note that $\sum_i d_i^2=\sum_j j^2n_j$, so Conditions \refand{C1}{C2}
imply, by Fatou's lemma, that $\sum_j j^2p_j<\infty$; in other words,
the asymptotic degree distribution has finite variance.

When Conditions \refand{C1}{C2} hold, the probability that \gndx{} is a
simple graph is bounded away from 0,
see \refSS{SSsimple} for details,
and thus all results that can be
stated in terms of convergence in probability for \gndx{} transfer to
the random simple graph \gnd{} too.

\subsection{Alternative graph models} 
We will in the remainder of the paper consider \gnd{}
as our underlying graph model, but we believe that similar results
hold for other random graph models too, and that they could be proved
by suitable modifications of the branching process arguments below.
Good candidates are the classical random graphs \gnp{} and \gnm{},
with $p=\mu/n$ and $m=n\mu/2$ (rounded to an integer), respectively,
and random graphs of the general type $\gnk$ defined in \cite{SJ178}.
We will not pursue this here, and leave such attempts to modify the
proofs to the interested reader,
but we will discuss one interesting case (including \gnp) where the
result easily follow from the results proved below for \gnd.

This example is a random graph defined by
Britton, Deijfen and Martin-L\"of \cite[Section 3]{BrittonDML},
see also \cite[Subsection 16.4]{SJ178}, as
follows.
Let $W$ be a non-negative random variable with finite expectation
$\mu_W\=\E W$. We first assign random weights $W_i$, $i=1,\dots,n$ to
the vertices; these weights are \iid{} with the same distribution as
$W$.
Secondly, given $\set{W_i}_1^n$, we draw an edge between vertices $i$
and $j$ with probability
\begin{equation}\label{pij}
  \pij\=\frac{W_iW_j}{n+W_iW_j};
\end{equation}
this is done independently (conditioned on $\set{W_i}$) for all pairs
\set{i,j} with $1\le i<j\le n$.
We denote this random graph by \gnw.
It is easily seen \cite{BrittonDML}
that \eqref{pij} implies that all graphs with a
given degree sequence $(d_i)_1^n$ have the same
probability. Hence, if we denote the (random) vertex degrees by
$D_1,\dots,D_n$, then conditioned on $D_i=d_i$, $i=1,\dots,n$, we have
a random graph \gnd.
Moreover, it is not difficult to verify that \refCN{C1} holds in
probability, with $(p_j)_0^\infty$ the mixed Poisson distribution
$\Po(\mu_WW)$ and $\mu=\mu_W^2$, see
\cite[Theorem 3.1]{BrittonDML} and
\cite[Theorem 3.13]{SJ178};
in other words,
$n_j/n\pto p_j$ and
$2m/n=n^{-1}\sum_id_i\pto\mu$.
Assume from now on that $\E W^2<\infty$; it may
then be shown by similar arguments that
$n^{-1}\sum_id_i^2\pto\mu_w^2(\E W^2+1)$.
Using the Skorohod coupling theorem,
see \eg{} \cite[Theorem 4.30]{Kallenberg}),
we can assume that these limits hold a.s.; hence Conditions
\refand{C1}{C2} hold a.s.
Consequently, by conditioning on $(D_1,\dots,D_n)$, we can apply the
results proved in the present paper for \gnd, and it follows that
the theorems below hold for the random graph \gnw{} too, with $(p_j)$
and $\mu$ as above.

With suitable couplings, using for example $(1+n^{-1/2})W_i$ for
upper bounds, it is easy to see that this remains true if \eqref{pij}
is modified to
\begin{equation}\label{pij1}
  \pij\=\min\Bigpar{\frac{W_iW_j}{n},1}.
\end{equation}
Random graphs defined by this definition and minor variations of it
have been studied by
several authors, see \cite[Subsection 16.4]{SJ178} and the references
given there.
Note that the special (deterministic) case $W=\sqrt\mu$ for a
constant $\mu>0$ gives the classical random graph $G(n,\mu/n)$.
The results in this paper thus holds for $G(n,\mu/n)$ too, with
$(p_j)$ a $\Po(\mu)$ distribution;
in other words, with $D$ defined in \refS{Sresults}, $D\sim\Po(\mu)$.

\subsection{Epidemic model}
We consider an infectious disease that spreads along the edges of
a graph $G$. We will in this paper assume that 
$G=\gnd$ is the random graph defined above, where we condition the graph
\gndx{} on being simple.
The vertices of $G$ are the individuals in the population, and the
edges represent friendships through which infection might spread.

The disease has its course in the following way.
Initially, one randomly chosen individual (vertex) is infected from the
outside. This individual then spreads the disease to each of its
friends independently and with the same probability $p$.
Those who get infected make out the first generation infected in
the epidemic. These individuals then do the same thing to their
not yet infected friends thus infecting a second generation,
and so forth. Note that an individual can only get infected once
--  we then consider such an individual either recovered and immune (or dead).
This epidemic continues until there are no new infections in a
generation, when it stops. Since the population is finite this happens after a
finite number of generations ($\le n$, where $n=|G|$ is the size of the
population). The individuals who get
infected during the course of the epidemic make up the total
outbreak, and the number of such individuals is called the final
size of the epidemic.

Note that each edge is a possible path of infection at most once,
namely when the first of its endpoints has been infected.
Hence we may just as well determine in advance for every edge in
$G$ whether it will spread the disease or not, provided that one of
the endpoints gets infected. Equivalently, we may consider the graph $G_p$
obtained by randomly deleting edges from $G$, with each edge kept with
probability $p$, independently of the others.
The final size of the epidemic is thus the size of the component of $G_p$
containing the initially infected individual.

\subsection{Vaccination strategies}
Assume now that a perfect vaccine is available. By this we mean
that an individual who is vaccinated is completely protected from
(\ie, immune to) the disease and is not able to spread the
disease further. We assume that a part of the population is vaccinated
before the epidemic starts, or as soon as the first individual is
infected. The epidemic progresses as defined above, with
the only difference that infected individuals can only infect
unvaccinated friends.

Note that for the study of the epidemic in the vaccinated population,
we may simply remove all
vaccinated individuals from $G$ (and edges connected to these
individuals). If we let $G\vv$ denote the remaining
graph,
and we assume that the initially infected individual $x$ is not vaccinated,
the final size of the epidemic is thus the size of the component of
$G\vvp\=(G\vv)_p$ that contains $x$.
We thus have to study the combined effect on $G$ of vertex deletion by the
vaccination and edge deletion by the randomness of infection.

The goal is to contain the disease, so that the final size of the
epidemic is small, and it is preferable to do this with a rather small
number of vaccinations.
For this we look at different local vaccination strategies.
The first two strategies are local in the sense that they require no global
knowledge of the social network $G$ (which is rarely available in
applications, \cite[Section 8.2]{n03}) and the latter two selects
friendhsips rather than individuals at random which may also be
thought of needing only local information.
We let $V$ denote the (usually random) number of vaccinations.

\subsubsection{Uniform vaccination}
Let us assume that we sample a fraction $c\in[0,1]$ chosen uniformly in the
population without replacement and that this fraction is
immunized, so the fraction $v$ being immunized satisfies $v=c$. This
vaccination strategy is the most commonly studied
vaccination strategy due to its simplicity  \cite[Section 8.2]{n03}.

More precisely, for convenience, we assume that each individual is
vaccinated with a given probability $v$, independently of each other.
The number $V$ of vaccinations is thus $\Bi(n,v)$, and $V/n\pto v$ as
\ntoo{} (with $v$ fixed).
We denote the remaining graph of unvaccinated individuals by $G\uv$;
this is thus obtained from $G$ by random vertex deletions.
Remember that our main concern is with the graph $G\uvp=(G\uv)_p$;
this is obtained from $G$ by random vertex and edge deletions,
independently for all vertices and edges. (In this case, it does not
matter whether we delete edges or vertices first.)

\subsubsection{Acquaintance vaccination}
It is intuitively clear that a better vaccination strategy would be to
vaccinate the individuals with
highest degrees (most friends) since this would reduce potential
spread the most.
However, for this targeted vaccination strategy to be achievable the
whole social graph (or at least the degrees of all individuals)
would have to be known, and this is rarely the case \cite[Section 8.2]{n03}.
A different
strategy aiming at vaccinating individuals with high degree, but
still only using local graph-knowledge from selected individuals,
proposed by Cohen et al.\ \cite{chb03},
goes under the name \emph{acquaintance vaccination}.
In this vaccination strategy a fraction $c$ of
individuals are sampled, and for each sampled individual one of
its friends, chosen randomly among all friends, is vaccinated. Of
course it may happen that some individuals are chosen more than
once for immunization (being selected as friends of more than one
individual) so the fraction $v=v(c)$ actually immunized is smaller than
$c$. This vaccination strategy has two slightly different variants
depending on whether the "fraction" $c$ is chosen with or without
replacement. We will use the version with replacement. For this case
the "fraction" $c$ may
in fact exceed 1 without having everyone vaccinated (individuals
who are selected more than once are asked for friends
independently each time and friends not yet immunized are
vaccinated). To be precise, we let the
number of individuals sampled be Poisson distributed $\Po(cn)$, with
$c\in[0,\infty)$.
Equivalently, each individual is sampled $\Po(c)$ times, and each time
reports a randomly chosen friend. Again, for simplicity, we assume
that each individual does this with replacement. Consequently, an
individual with degree $d$ will report each of its friends
$\Po(c/d)$ times, and these random numbers are all independent.
(An individual that is sampled but has no friends is ignored.
An individual is only vaccinated once, even if he or she is reported
several times.)

For any initial graph $G$ and $0\le c<\infty$, 
we denote the remaining graph of unvaccinated individuals by $G\ac$.
We further write  $G\acp=(G\ac)_p$ for the graph obtained by additional
edge deletions. (For acquaintance vaccination, the order of the
deletions is important, since the vaccination strategy uses all edges,
without knowing whether they may be selected to transmit the disease
or not.)

\subsubsection{Edgewise vaccination} 
In some situations it may be possible to observe, or at least sample,
the edges representing friendships. If this is the case, another
reasonable vaccination strategy is to sample a number of the edges and
then either vaccinate both endpoints or one (randomly selected)
endpoint; we denote these two versions by $\eix$ and $\eiix$.

For $\eiix$, we assume that we sample each edge with probability
$1-\alpha$, where $\alpha\in(0,1]$ is a fixed number. (Equivalently,
we sample $\Po(cm)$ edges with replacement, with $\alpha =e^{-c}$.)
For $\eix$, we assume for simplicity that we sample $\Po(2cm)$ edges
with replacement; thus each end of each edge is sampled with
probability $1-\alpha =1-e^{-c}$, independently of all other edge
ends. Hence, for both versions, a vertex with degree $d$ is
unvaccinated with probability $\alpha^d$, and for $\eix$, this is
independent of all other vertices.

For an initial graph $G$ and $0< \alpha\le1$,
we denote the remaining graph of unvaccinated individuals by $G\eiq$
and $G\eiiq$, for the two versions.
We further write, for $j=1,2$,
$G\ejqp=(G\ejq)_p$ for the graph obtained by additional
edge deletions.

\section{Main results}\label{Sresults}

We now state our main results together with heuristic
motivations.
We assume that the underlying graph is the random graph $\gnd$ and
that Conditions \refand{C1}{C2} hold.
Complete proofs are given in Section \ref{Spf}.

\subsection{Original epidemic model}\label{SSomodel}

Assume that $n$, the number of nodes, is large. The regularity
assumption 
on the degrees of
the graph (\refCN{C1}) implies that no
separate node will contain a large fraction of all edges, see
\eqref{unif1}.
This in turn
implies that self loops, multiple edges and short cycles will be rare.

The epidemic starts by a randomly
selected individual being infected from outside, so this individual
has (approximately) 
the degree
distribution $(p_j)_{j=0}^\infty$. The friends of this individual,
or of any individual, have the size
biased degree distribution $(\tilde p_j)_{j=0}^\infty$, where
\begin{equation}
\tilde p_j=jp_j/\sum_kkp_k.\label{ptilde}
\end{equation}
Let $D$ and $\tD$ be random
variables having these degree distributions respectively. Then, given that
$D=d$, the number of
individuals that the initially infected infects is $\Bi (d,p)$, and
the unconditional distribution is hence mixed binomial $\MixBi
(D,p)$. Those then infected, as well as infecteds
in the following
generations, have degree distribution $(\tilde
p_j)_{j=0}^\infty$. Given that $\tD=\tilde d$,
the number of individuals an infected individual infects in the next
generation has distribution $\Bi(\tilde d-1,p)$. This follows because
the infected was infected by one of his friends (which cannot get
reinfected) and, since short cycles are
rare, it is very unlikely that any of the remaining $\tilde d-1$
friends have already been infected. Unconditionally, the number
infected in the next generation is hence $\MixBi (\tD-1,
p)$. Further, the property that short
cycles are unlikely implies that the number of infections caused by
different individuals are (approximately) independent random variables.

The above paragraph motivates why the early stages of the epidemic
may be approximated by a branching process (e.g.\ \cite{AN}), as is
common for epidemic models (e.g.\ \cite{ab00}), and where ``giving birth''
corresponds to infecting someone. The branching process is a simple
Galton--Watson process
starting with one ancestor having off-spring distribution  $X\sim
\MixBi(D,p)$ and the following generations have off-spring distribution
$\tilde X\sim \MixBi(\tD-1,p)$. The
mean of this latter off-spring distribution plays an important role
in branching process theory and also in in epidemic theory where it is
denoted $R_0$ and denoted the basic reproduction number.  We
get the following, using \eqref{ptilde}, 
\begin{equation}
  \begin{split}
R_0=\E(\tilde X)=p\E(\tD-1)
=p\lrpar{\frac{\sum_jj^2p_j}{\mu}-1}
=p\left(\mu+\frac{\Var(D)-\mu}{\mu}\right),\label{R_0}
  \end{split}
\end{equation}
 where $\mu=\E(D)=\sum_kkp_k$ and $\Var(D)=\sum_jj^2p_j - \mu^2$ (a
very related expression is obtained in \cite{a99}).
The branching process is subcritical, critical or supercritical
depending on whether $R_0<1$, $R_0=1$ or $R_0>1$. For the epidemic,
this means a major outbreak infecting a non-negligible fraction of the
community, is possible if and only if $R_0>1$. Note that, for fixed
$\mu$, $R_0$ is increasing in $\Var(D)$, so the more variance in the
degree distribution, the higher $R_0$, and if the degree distribution
has infinite variance then $R_0=\infty$ (a case not treated in the
present manuscript due to Condition \ref{C2}).

The probability $\pi$ that the branching process dies out is derived
in the standard way 
as follows.
First, we derive the probability $\tilde \pi$ that a branching process
with all individuals having off-spring distribution $\tilde X$ dies
out. This is obtained by conditioning on the number of individuals
born in the first generation: for the branching process to die out,
all branching processes initiated by the individuals of the first
generation must die out, i.e.\
$$
\tilde \pi=\sum_{k=0}^\infty \tilde \pi^k\P(\tilde X=k).
$$
Let $f_{\tilde X}(\cdot)$ denote the probability generating function for
$\tilde X$,  and $f_{ D}(\cdot)$ the probability generating function
of the original degree distribution $D$. Then we see
that $\tilde \pi$ is a solution to the equation $f_{\tilde X}(t)=t$,
and it is known from branching process theory (e.g. \cite[Theorem
I.5.1]{AN}) that it is
the smallest non-negative 
such solution. The fact that $\tilde X$ is $\MixBi(\tilde
D-1,p)$ implies that
$$
f_{\tilde X}(t)=\E(t^{\tilde X})=\E(\E(t^{\tilde X}|\tilde
D))=\E((pt+1-p)^{\tD-1})=\E((1-p(1-t))^{\tD-1}).
$$
Further,
$$
\E(a^{\tD-1})=\sum_ka^{k-1}\frac{kp_k}{\mu}=
\frac{d}{da}\sum_ka^k\frac{p_k}{\mu}=\frac{d}{da} \frac{f_D(a)}{\mu }
=\frac{f_D'(a)}{\mu}=\frac{f_D'(a)}{f_D'(1)}.
$$
In terms of $f_{ D}(\cdot)$ the probability $\tilde \pi$ that the branching
process dies out is hence the smallest non-negative solution to
\begin{equation}
\frac{f_D'(1-p(1-\tilde \pi))}{f_D'(1)}=\tilde \pi.\label{tildepi}
\end{equation}

The probability $\pi$ that the branching process, in which the
ancestor has different off-spring distribution $X$, dies out, is obtained from
$\tilde \pi$ by conditioning on the number of off-spring of the
ancestor:
\begin{equation}
  \begin{split}
\pi&=\sum_k\tilde \pi^k\P(X=k)=\E(\tilde \pi^X)=\E(\E(\tilde
\pi^X|D))=\E((p\tilde \pi+1-p)^D)
\\&
=f_D(1-p(1-\tilde \pi)).\label{pi}
  \end{split}
\end{equation}

We now look at the final size of the epidemic in case it takes off,
corresponding to the case that the branching process grows beyond
all limits. We do this by considering the epidemic from a graph
representation. The social structure was represented by a random
graph $G$. If this graph is thinned by removing each edge
independently with probability $1-p$ we get a thinned graph denoted
$G_p$. Edges in $G_p$ represent potential spread of infection: if
one of the nodes get infected from elsewhere, its neighbour will get
infected. As a consequence, the final outbreak of the epidemic will
consist of all nodes in $G_p$ that are connected to the initially
infected. From random graph theory it is known that if $R_0>1$ there
will be exactly one connected component of order $n$, the giant
component, and all remaining connected components will be of smaller
order. If $R_0\le 1$ there will be no giant component. The initially
infected was chosen uniformly in the community so it will belong to
the giant component with a probability that equals the relative size
of the giant component. On the other hand, the initially infected
belongs to the giant component if and only if its branching process
of new infections grows beyond all limits, and we know from before
that this happens with probability $1- \pi$ defined in equation
(\ref{pi}). From this it follows that the asymptotic final
proportion infected, $\tau$, equals $1- \pi$. So, $\tau$ is both the
probability of a major outbreak, and the relative size of the
outbreak in case a major outbreak occurs.

The above arguments motivate the following theorem, which is proven in
Section \ref{Spf},
and where $Z_n$ denotes the final number infected in the epidemic.
\begin{theorem}\label{T0}
If $R_0\le 1$ then $Z_n/n\pto0$. If $R_0>1$, then $Z_n/n$ converges to a
two-point distribution $Z$ for which $\P(Z=0)=\pi$ and $\P(Z=\tau)=\tau$,
where $\pi$ is defined by \eqref{tildepi} and \eqref{pi} and $\tau =1-\pi$.
\end{theorem}

\subsection{Uniform vaccination} \label{SSumodel}
Prior to arrival of the infectious disease, each individual is
vaccinated independently and with the same probability $v$ which
implies that the total number of vaccinated $V$ is $\Bi(n,v)$, and
from the law of large number the random proportion vaccinated
$V/n\pto v$.

Vaccinated individuals, and edges connecting to them, can be removed from the
graph since there will be no spreading between these
individuals and their friends in either direction. As a consequence,
an individual who originally had $d$ friends now has $\Bi(d,1-v)$
unvaccinated friends. If an individual gets infected during the
early stages of the epidemic he will infect
each of his unvaccinated friends independently with probability
$p$. Given that the initially infected has degree
$d$ he will hence infect $\Bi(d,p(1-v))$ friends, so without the
conditioning he will infect a mixed binomial number $X_v\sim
\MixBi(D,p(1-v))$. Similarly, during the
early stages an infected individual with degree $d$
will infect $\Bi(d-1,p(1-v))$, and unconditionally an individual has
degree distribution $\{\tilde p_k\}$, so the unconditional number he
will infect $\tilde X_v$ will be $ \MixBi(\tilde
D-1, p(1-v))$.

It is seen that we have the same type of distributions as in the case
without vaccination. As a
consequence, all results for the case with uniform vaccination can be
obtained from the case without vaccination simply by replacing $p$ by
$p(1-v)$. We hence have that the reproduction number $R\uvp$ after
vaccinating a
fraction $v$ chosen uniformly satisfies
\begin{equation}
R\uvp=\E(\tX_v)=(1-v)R_0
=p(1-v)\left(\mu+\frac{\Var(D)-\mu }{\mu}\right).\label{R_vU}
\end{equation}
The probability $\tilde \pi\uvp$ that the epidemic never takes
off, assuming the
initially infected has $\tilde X_v$ unvaccinated friends, is
the smallest
solution to
\begin{equation}
\frac{f_D'\bigpar{1-p(1-v)(1-\tilde \pi\uvp)}}{f_D'(1)}=\tilde \pi\uvp.\label{tildepi_vU}
\end{equation}
The probability $\pi\uvp$ that the epidemic never takes off if
the initially infected is selected randomly among the unvaccinated is
given by
\begin{equation}
\pi\uvp=f_D\bigpar{1-p(1-v)(1-\tilde \pi\uvp))},\label{pi_vU}
\end{equation}
where $\tilde \pi\uvp$ is the smallest solution to (\ref{tildepi_vU}).
Finally, the final size is determined from the probability of a major
outbreak as before. This means that the final proportion infected
(among the unvaccinated!)\ will converge to $1-\pi\uvp$ in case of a
major outbreak. We have the following corollary, where $Z_n\uq(v)$
denotes the final number infected in the epidemic where each
individual was vaccinated independently with probability $v$ ($0\le
v<1$) prior to
the outbreak, and where the initially infected was chosen randomly
among the unvaccinated.
\begin{theorem}\label{TU}
If\/ $R\uvp\le 1$, then $Z_n\uq(v)/((1-v)n)\pto0$. 
If\/ $R\uvp>1$,
then $Z_n\uq(v)/((1-v)n)$ 
converges to a
two-point distribution $Z\uvp$ for which
$\P(Z\uvp=0)=\pi\uvp$ and\/
$\P(Z\uvp=\tau\uvp)=\tau\uvp$,
where $\pi\uvp$ is defined by \eqref{tildepi_vU} and
\eqref{pi_vU} and $\tau\uvp =1-\pi\uvp$.
\end{theorem}

\subsection{Acquaintance vaccination}\label{SSamodel}
Recall that each individual is sampled, independently, a $\Po(c)$
number of times, where $0\le c<\infty$, so in total $\Po(nc)$
individuals are sampled. Each time an individual is sampled a
randomly chosen friend of the individual is selected and vaccinated
(unless it already was vaccinated). The effect of this strategy is
that vaccinated individuals have the size biased degree distribution
$(\tilde p_j)_{j=0}^\infty$, where $\tilde p_j=jp_j/\sum_kkp_k$
rather than the original degree distribution $\{ p_k\}$ for
uniformly selected individuals. The proportion vaccinated $v=v(c)$
is obtained as follows. An individual avoids being vaccinated if he
is not vaccinated ``through'' any of its friends. The friends of the
individual have independent degree distributions $(\tilde
p_j)_{j=0}^\infty$, and the probability of not being vaccinated
``through'' an individual with degree $k$ is $e^{-c/k}$. It follows
that the probability to avoid being vaccinated from one friend
equals
\begin{equation}
\ga=\ga(c)=\sum_{k=1}^\infty e^{-c/k}\tilde p_k =\sum_{k=1}^\infty
e^{-c/k}\frac{kp_k}{\mu}.\label{alpha}
\end{equation}
(Note that $\alpha$ has the same interpretation as for $\alpha$
introduced for the edgewise strategies, but it is a different
function of $c$.) If the individual in question has $j$ friends it
hence avoids being vaccinated with probability $\ga^j$. The
proportion $1-v(c)$ not being vaccinated equals the probability that
a randomly selected individual is \emph{not} vaccinated, which hence
equals
\begin{equation}\label{v_A}
1-v(c)=\sum_{j=0}^\infty \ga^jp_j=f_D(\ga),
\end{equation}
where as before $f_D(\cdot )$ is the probability generating function
of a random
variable $D$ having distribution $(p_j)_{j=0}^\infty$.

Note that in this model, given the graph, individuals are vaccinated 
independently of each other (although with different probabilities).
It follows easily that the actual (random) number $V$ of vaccinated
persons satisfy
\begin{equation}\label{V_A}
  V/n \pto v(c)
\qquad \text{as \ntoo}.
\end{equation}
Hence we will ignore the randomness in $V$ and regard $v(c)$ given by
\eqref{v_A} as the proportion of vaccinated persons.

We now approximate the initial stages of an epidemic, occurring in a
community having been
vaccinated according to the acquaintance strategy, with a suitable
branching process. To find ``the right'' branching process
approximation is harder
for the acquaintance strategy because the vaccination status of an
individual depends on the
degrees of its friends. We therefore introduce some convenient
terminology.

We say that \emph{transmission may take place} through an edge, and
 through its two half-edges, if it is one of the edges in $G_p$,
i.e., one of the randomly selected edges which will spread the
disease if one of its endpoints is infected. (Recall that we may
assume this random selection to take place before the start of the
infection.) Further, there is a natural correspondence between
half-edges and \emph{directed} edges, with a half-edge corresponding
to the edge it is part of, directed so that the it begins with this
half-edge. We say that a directed edge, or the corresponding
half-edge, is \emph{used for vaccination}, if the person at the
start of the edge is selected and names the person at the end of the
edge, who thus gets vaccinated.

\begin{figure}[htb] \begin{center} \bf
\includegraphics[width=12.0cm, height=4.5cm]{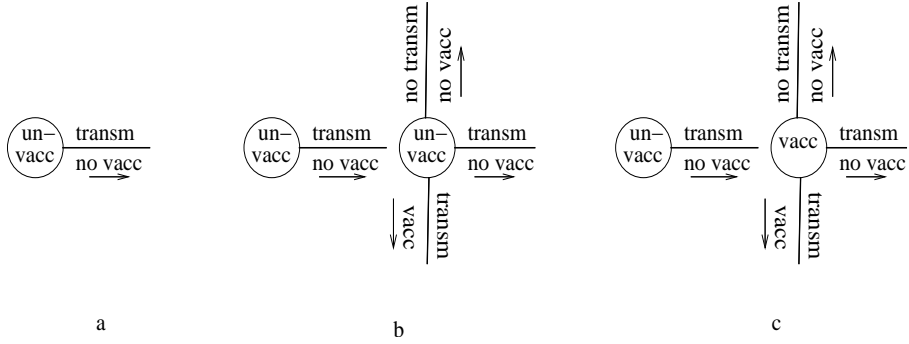}
 \caption{\rm a) An illustration of an ``individual'' in the branching
 process. In b) the left ``individual'' has one off-spring (the up-going
edge constitutes no individual since there is no transmission and the
down-going no individual since the friend was sampled and
named the individual below for vaccintation. In c) no individual is born
since the friend was vaccinated (being named by some other friend).}
\label{acq-br}\end{center}
\end{figure}

It turns out that a suitable ``individual'' in the branching process
is an unvaccinated person \emph{together} with a directed edge from
this person such that transmission may take place through the edge
but it is not used for vaccination. It is worth noting that a person
may be part of several ``individuals'' in the branching process (if
the person was not vaccinated and has several friends such that the
connecting edges satisfy the conditions above).  See Figure
\ref{acq-br} for an illustration of an individual (a) and situations
where the individual ``gives birth'' to one (b) and 0 (c)
individuals.

In order to analyse the corresponding branching process we have to
determine the distribution of how many new
``individuals'' one ``individual'' will infect during the early stages
of the epidemic assuming a large population (large $n$). We know that
the ``individual'' contains an unvaccinated person, so
the edge in the ``individual'' has
not been used for
vaccination backwards, \ie{} in
the opposite direction. 
As a consequence, we have to condition on this, and then the node
at the other end of the edge
has degree $K=k$
with probability
\begin{equation}
  \label{pfriend}
\P(K=k)=\frac{\tilde p_ke^{-c/k}}{\sum_{j=1}^\infty \tilde
  p_je^{-c/j}}=\frac{\tilde p_ke^{-c/k}}{\ga }, \quad k=1,2,\ldots,
\end{equation}
i.e.\ the size biased degree distribution conditional on not having
vaccinated backwards. In order for this friend to create new
``individuals'', it must not have been vaccinated by any of his other
$k-1$ friends (by assumption it was not vaccinated from our original
individual). This happens with probability $\ga^{k-1}$. Each of the friend's
remaining $k-1$ edges will be \emph{open} (i.e., transmission may
take place but it is not used for vaccination) 
\emph{independently}, each open with probability $pe^{-c/k}$. The number
of open edges (equal to the number of new ``individuals'') is
hence $\Bi(k-1,pe^{-c/k})$. If the
friend is vaccinated (probability $1-\ga^{k-1}$) no new individuals are
born.  The unconditional number $Y$ of new
``individuals'' an individual ``gives birth'' to, i.e.\ the off-spring
distribution of the approximating branching process, can be obtained
by conditioning on the number of friends our friend has and recalling
that 0 individuals are born whenever the friend is vaccinated or 
if the binomial variable equals 0:
\begin{equation}
  \label{aoff}
\begin{aligned}
\P(Y=0)&= \sum_{k= 1}^\infty \left( (1-\ga^{k-1})
  + \ga^{k-1}(1-pe^{-c/k})^{k-1}\right)\frac{\tilde p_ke^{-c/k}}{\ga} ,
\\
\P(Y=j)
&= \sum_{k= j+1}^\infty \ga^{k-1}\binom{k-1}{j}
 (pe^{-c/k})^j(1-pe^{-c/k})^{k-1-j}\frac{\tilde p_ke^{-c/k}}{\ga},\quad j\ge 1.
\end{aligned}
\end{equation}
 This
off-spring distribution determines both $R\acp$,
the probability of a
major outbreak, and the final size in case of a major outbreak. For
instance, the reproduction number is the mean of this distribution,
and this mean is obtained by first conditioning on the degree of the
node in question. Given that the degree equals $k$, the average number of
off-spring equals $\ga^{k-1}(k-1)pe^{-c/k}$, which gives the following
reproduction number:
\begin{equation}
R\acp
=\E(Y)
=
\sum_{k\ge 1} \ga^{k-1}(k-1)pe^{-c/k}\frac{\tilde
  p_ke^{-c/k}}{\ga} =p\sum_{k\ge 1}(k-1) \ga^{k-2}e^{-2c/k}\tilde p_k
\label{R_vA}
\end{equation}
(cf.\ \cite{chb03}). Let $f_Y (a)=\E(a^Y)$ be the probability
generating function of this off-spring distribution. If the epidemic
starts by one ``individual'', \ie{} one person with one open
directed edge,  then the probability $\tilde \pi\acp$ that the
epidemic never takes off is the smallest solution to the equation
\begin{equation}
\tilde \pi\acp=f_Y (\tilde \pi\acp).\label{tildepi_vA}
\end{equation}
If we start with one infected person that is unvaccinated and has
 degree $j$, then each of its $j$ half-edges is open with
probability $pe^{-c/j}$, and the probability that a given half-edge
does not start a large epidemic is
$1-pe^{-c/j}+pe^{-c/j}\tilde\pi\acp$, so the probability that the
epidemic never takes off equals $(1-pe^{-c/j}(1-\tilde\pi\acp))^j$,
for $j\ge 1$, and 1 for $j=0$.

If instead the initially infected is chosen randomly
among the unvaccinated as we assume,
then the probability that it has degree $j$ is
$p_j \ga^j/\sum_j p_j \ga^j$,
\cf{} \eqref{v_A}, and thus the probability that the
epidemic never takes off equals
\begin{equation} \label{pi_vA}
\pi\acp=
\frac{p_0+\sum_{j\ge 1} p_j \ga^j\bigpar{1-pe^{-c/j}(1-\tilde \pi\acp)}^j}
{\sum_j p_j \ga^j}
.
\end{equation}
Finally, using the same reasoning as before, the limiting proportion
infected in case of a major outbreak equals $\tau\acp=1- \pi\acp$.
We summarize our results in the following theorem, proved in
\refS{Spf}, where $Z_n\aq\vc$ denotes the final number infected in
the epidemic where vaccination is done prior to the outbreak
according to the acquaintance vaccination strategy. Recall that
$0\le c<\infty$ and that $v(c)$, the proportion of the population
vaccinated, is given by \eqref{v_A} with $\ga=\ga(c)$ given by
\eqref{alpha}.
\begin{theorem}\label{TA}
 $Z_n\aq\vc/\bigpar{(1-v(c))n}\pto0$  
if
$R\acp\le 1$, where $R\acp$
is defined by \eqref{R_vA}. If $R\acp>1$, then
$Z_n\aq\vc/\bigpar{(1-v(c))n}$ 
converges to a
two-point distribution $Z\acp$ for which
$\P\bigpar{Z\acp=0}=\pi\acp$ and
$\P\bigpar{Z\acp=\tau\acp}=\tau\acp$,
where $\pi\acp$ is defined by \eqref{tildepi_vA} and
\eqref{pi_vA}, and $\tau\acp =1-\pi\acp$.
\end{theorem}

\subsection{Edgewise vaccination}\label{SSemodel}
Recall that, for both $\eix$ and $\eiix$, a person with $d$ friends
is unvaccinated with probability $\ga^d$ (here $\ga$ has the same
meaning in the previous subsection, but it can be treated as a free
parameter). Thus,
\begin{equation*}
\E V  = n\sum_d p_d(1-\ga^d) +o(n)
\end{equation*}
and a simple variance estimate shows that the vaccinated proportion
\begin{equation}\label{v_E}
  V/n \pto v(\ga):=\sum_d p_d(1-\ga^d),
\end{equation}
just as for acquaintance vaccination, see \eqref{v_A} and \eqref{V_A}.

We define open (directed) edges as for acquaintance vaccination, and
argue as there with the following modifications. The other endpoint of
an open edge has just the size-biased distribution $(\tp_k)$. If this
vertex, $y$ say, has degree $k$, it is unvaccinated with probability
$\ga^{k-1}$, and in that case, the number of new open edges originating
at $y$ is $\Bi(k-1,p\ga)$ for $\eix$ and $\Bi(k-1,p)$ for $\eiix$.
The difference between the two versions is because we already know
that these edges do not vaccinate $y$, and for $\eiix$, this implies
that they do not vaccinate their other endpoint either, while for
$\eix$ that is an independent event with probability $\ga$.

We thus have the offspring distributions for $\eix$ and $\eiix$,
\cf{} \eqref{aoff},
\begin{equation*}
\begin{aligned}
\P(Y_1=j)
&= \sum_{k= j+1}^\infty \tp_k \ga^{k-1}\binom{k-1}{j}
 (p\ga)^j(1-p\ga)^{k-1-j}, &&\quad j\ge 1,
\\
\P(Y_2=j)
&= \sum_{k= j+1}^\infty \tp_k \ga^{k-1}\binom{k-1}{j}
 p^j(1-p)^{k-1-j}, && \quad j\ge 1;
\end{aligned}
\end{equation*}
we leave the formulas for $\P(Y_1=0)$ and $\P(Y_2=0)$ to the reader.

This gives the reproduction numbers
\begin{equation}\label{R_vE}
  \begin{aligned}
R\eiqp
&
=\E(Y_1)
=
\sum_{k\ge 1}\tp_k \ga^{k-1}(k-1)p\ga
=p\sum_{k}(k-1)\tp_k \ga^{k},
\\
R\eiiqp
&
=\E(Y_2)
=
\sum_{k\ge 1}\tp_k \ga^{k-1}(k-1)p
=p\sum_{k}(k-1)\tp_k \ga^{k-1}.
  \end{aligned}
\end{equation}
Note that $R\eiqp=\ga R\eiiqp<R\eiiqp$, which shows that, with the same
number of vaccinations, $\eix$ is a better strategy than $\eiix$. 
In particular, the critical
critical vaccination coverage $v_c$ is smaller for 
$\eix$ than for $\eiix$. 
An intuitive 
explanation to why $\eiix$ is not as efficient as $\eix$ is that in
$\eiix$ both individuals of selected friendships are vaccinated, and
since an individual is partly protected by friends getting
vaccinated the second vaccination is less ``efficient''.

We let $\tpi\eiqp$ and $\tpi\eiiqp$ be the probabilities that the
Galton--Watson processes with offspring distributions $Y_1$ and
$Y_2$, respectively, starting with one individual, die out; they are
thus the smallest positive solutions to $t=f_{Y_1}(t)$ and
$t=f_{Y_2}(t)$, where $f_{Y_1}$ and $f_{Y_2}$ are the corresponding
probability generating functions.

If we start with one unvaccinated person $x$ with degree $d$,
the number of open edges from $x$ is $\Bi(d,p\ga)$ for $\eix$ and
$\Bi(d,p)$ for $\eiix$, for the same reason as for the number of new
edges above.
The probability that the
epidemic never takes off is thus
$(1-p\ga+p\ga\tpi\eiqp)^d$ for $\eix$ and
$(1-p+p\tpi\eiiqp)^d$ for $\eiix$.

If the initially infected is chosen randomly
among the unvaccinated,
we thus find the probabilities that the
epidemic never takes off
\begin{equation} \label{pi_vE}
  \begin{aligned}
\pi\eiqp &= \frac{\sum_j p_j \ga^j\bigpar{1-p\ga(1-\tpi\eiqp)}^j}
{\sum_j p_j \ga^j} =
\frac{f_D\bigpar{\ga\bigpar{1-p\ga(1-\tpi\eiqp)}}} {f_D(\ga)},
\\
\pi\eiiqp &= \frac{\sum_j p_j \ga^j\bigpar{1-p(1-\tpi\eiiqp)}^j}
{\sum_j p_j \ga^j} =
\frac{f_D\bigpar{\ga\bigpar{1-p(1-\tpi\eiiqp)}}} {f_D(\ga)} .
  \end{aligned}
\end{equation}

We summarize our results as before, letting
$Z_n\ei(\ga)$  and $Z_n\eii(\ga)$
denote the final numbers infected in the epidemic for the two
strategies. Recall that $v(\ga)$
is given by \eqref{v_E}.
\begin{theorem}\label{TE}
For $j=1,2$,
 $Z_n\ej(\ga)/\bigpar{(1-v(\ga))n}\pto0$
if
$R\ejqp\le 1$, where $R\ejqp$
is defined by \eqref{R_vE}. If $R\ejqp>1$, then
$Z_n\ej(\ga)/\bigpar{(1-v(\ga))n}$
converges to a
two-point distribution $Z\ejqp$ for which
$\P\bigpar{Z\ejqp=0}=\pi\ejqp$ and
$\P\bigpar{Z\ejqp=\tau\ejqp}=\tau\ejqp$,
where $\pi\ejqp$ is defined by
\eqref{pi_vE}, and $\tau\ejqp =1-\pi\ejqp$.
\end{theorem}

\subsection{Examples}\label{Ex}
\begin{figure}[htb] \begin{center} \bf
\includegraphics[width=13.0cm, height=10cm]{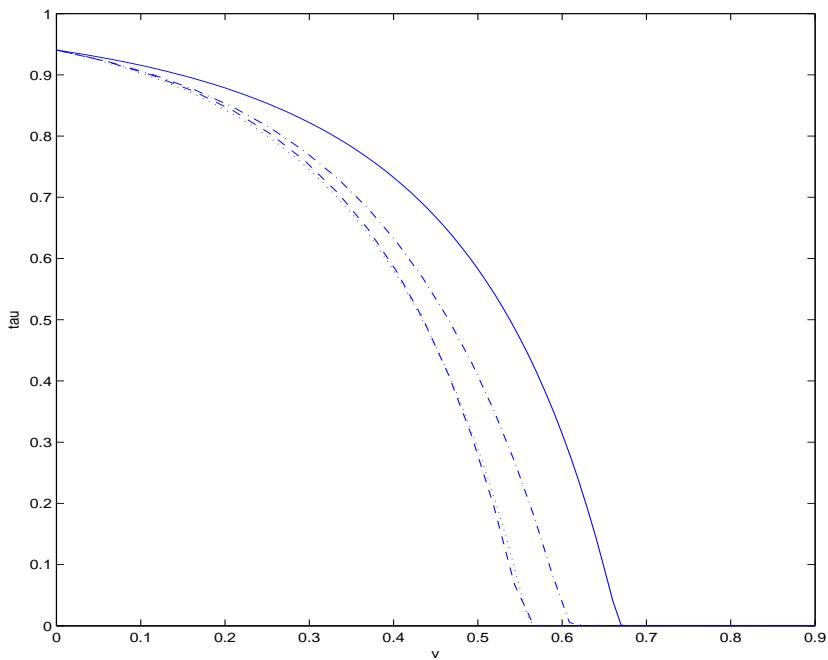}
 \caption{\rm Final proportion infected $\tau$ as a function of the
   vaccination coverage $v$ for four vaccination strategies: uniform
   (---), acquaintance ($\cdots$), $\eix$ (- - -) and
   $\eiix$ ($-\cdot -\cdot -$). The degree distribution is $\Po(6)$
 and transmission probability $p=0.5$.}
\label{tau-Po}\end{center}
\end{figure}

We now compare the performance of the different vaccination
strategies on two examples. In the first example we have chosen the
degree distribution to be Poisson distributed with mean $\lambda
=6$, and the transmission probability to equal $p=0.5$. Using
\eqref{R_0} we conclude that this implies that $R_0=3$. The
assumption of Poisson distributed degree means that this applies to
the 
simple $G(n,p=6/n)$ graph with transmission probability $p=0.5$; in
the epidemic literature this model is knowns as the Reed-Frost model
(e.g.\ \cite{ab00}). In Figure \ref{tau-Po} we show $\tau$, the
final proportion infected among unvaccinated in case of a major
outbreak, as a function of the vaccination coverage $v$, for the 4
different vaccination strategies treated.
It is seen that the acquaintance and edgewise $\eix$ strategies
perform best in the sense that, for a fixed proportion vaccinated,
the proportion $\tau$ getting infected in case of a major outbreak
is smallest for these two strategies. As a consequence, the critical
vaccination coverage, $v_c=\inf_v \{v; R_v\le 1\}$, is also smallest
for these two strategies. There is no unique ordering of the two
strategies -- the acquaintance strategy is slightly better for small
vaccination coverages and $\eix$ is slightly better for higher
vaccination coverages and hence also has slightly smaller $v_c$. The
edgewise strategy $\eiix$ is not as good as these two strategies but
still better than the uniform vaccination coverage. 
(Indeed, $\eiix$  is always less efficient than $\eix$, see above.) 
Acquaintance, $\eix$ and $\eiix$ all perform better than the uniform
strategy, the reason being that they tend to find individuals with
high degrees. For the parameter choices of this example, the
critical vaccination coverages equal $v_c\approx 0.56$ for the
acquaintence and  $\eix$ strategies, $v_c\approx 0.61$ for $\eiix$
and $v_c\approx 0.67$ for the uniform vaccination strategy.

\begin{figure}[htb] \begin{center} \bf
\includegraphics[width=13.0cm, height=10cm]{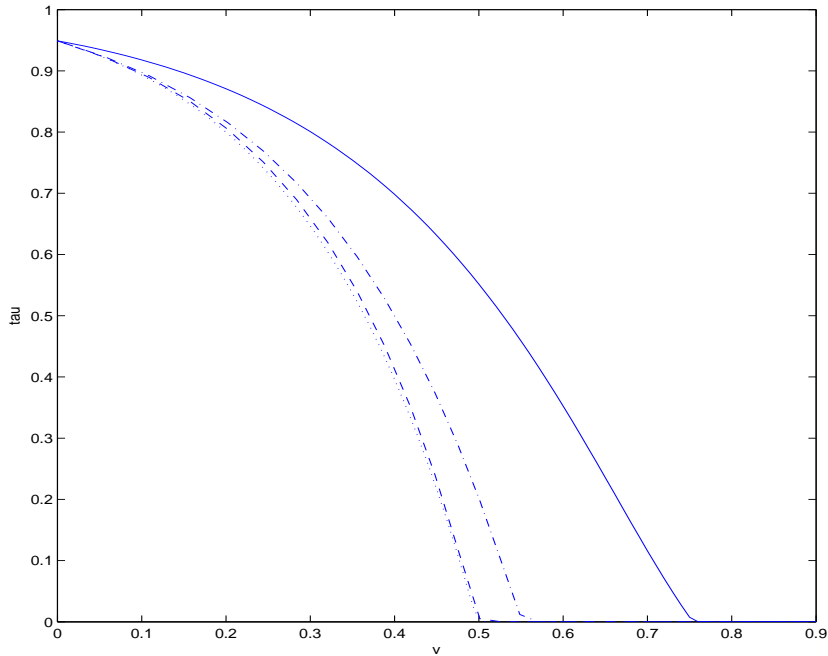}
 \caption{\rm Final proportion infected as a function of the
   vaccination coverage for four vaccination strategies: uniform
   (---), acquaintance ($\cdots$), $\eix$ (- - -) and
   $\eiix$ ($-\cdot -\cdot -$.) The degree distribution is
   heavy-tailed ($p_d\propto d^{-3.5}$) with mean $E(D)\approx 6$
 and $p=0.5$.}
\label{tau-heavy}\end{center}
\end{figure}

In the second example (illustrated in Figure \ref{tau-heavy}) we
chose a more heavy tailed degree distribution having $p_d\propto
d^{-3.5}$ (in the computations it was truncated at $d=200$). The
initial values were modified such that $E(D)\approx 6$ to make it
more comparable to the previous example, with a resulting variance
equal to 18.9. The transmission parameters was set $p=0.5$ as
before. Using \eqref{R_0} we hence see that $R_0\approx 4.1$.
In the figure we see the same type of pattern as in the previous
example. However, the difference between the strategies is more
pronounced with $v_c\approx 0.50$ for the acquaintence and $\eix$
strategies, $v_c\approx 0.55$ for $\eiix$ and $v_c\approx 0.75$ for
the uniform vaccination strategy. In other words, if the uniform
strategy is applied in these two examples we have to vaccinate
\emph{more} individuals if the degree distribution is heavy-tailed,
but if any of the other strategies is performed, the heavy-tailed
degree distribution require \emph{less} vaccinations to surely
prevent an outbreak. Another minor difference from the previous
example is that, for the present heavy-tailed distribution, the
acquaintance strategy is (slightly) better than $\eix$ for all
vaccination coverages and hence also has a smaller critical
vaccination coverage. However, the difference between the two
strategies is negligible.

Note that all $\tau$'s in both examples denote the proportion of
infected among the unvaccinated (in case of an outbreak) and can
hence be thought of as an indirect protection from those getting
vaccinated. Of course, by assumption, all vaccinated are also
protected from getting infected.

\section{Preliminaries on branching processes}
\label{Sgw}

As said above, our method is based on comparison with branching
processes, more precisely Galton--Watson processes, see \eg{} \cite{AN}
for definitions and basic facts.
If $\bp$ is a Galton--Watson process started with 1 initial particle,
we let $\bpd$ denote the same branching process with $d$ initial
particles, \ie{} the union of $d$ independent copies of $\bp$.
Further, for any Galton--Watson process $\bp$, we let
$|\bp|$ denote its total progeny, \ie{} the total number of particles
in all generations,
and we let $\rho(\bp)$ be the survival probability of $\bp$, \ie{}
$\rho(\bp)\=\P(|\bp|=\infty)$.
Note that if $\bp$ starts with 1 particle, then
\begin{equation}\label{erika}
  \rho(\bpd)=1-\bigpar{1-\rho(\bp)}^d,
\end{equation}
since $\bpd$ dies out if and only if all $d$ copies of $\bp$ in it do.

We will need the following simple continuity result, which presumably
is well known although we have failed to find a reference.

\begin{lemma}
  \label{Lgw}
Let $X_\nu$ and $X$ be non-negative integer-valued random variables,
and let $\bpnud$ and $\bpd$ be the corresponding
Galton--Watson processes with offspring distributions $X_\nu$ and $X$,
starting with $d$ particles.
If $X_\nu\dto X$ as \ntoo, and $\P(X=1)<1$, then
$\rho(\bpnud)\to\rho(\bpd)$,
for every fixed $d\ge0$.
\end{lemma}

\begin{proof}
By \eqref{erika}, it suffices to show this for $d=1$, and we then drop
the superscript 1.

Consider the probability generating functions
$f_X(t)\=\E t^X$ and $f_{X_\nu}(t)\=\E t^{X_\nu}$ for $0\le t\le1$.
It is well-known, see \eg{} \cite[Theorem I.5.1]{AN}, that
the extinction probability $q\=1-\rho(\bp)$ is the smallest root in $[0,1]$
of $f_X(q)=q$.
It follows easily, since we have excluded the possibility
$f_X(t)\equiv t$, that if $0\le t<q$, then $f_X(t)>t$,
and if $q<t<1 $, then $f_X(t)<t$.

Since $X_\nu\dto X$, we have $f_{X_\nu}(t)\to f_X(t)$ for every $t\in[0,1]$.
Hence, if $0\le t<q$, then $f_{X_\nu}(t)>t$ for large $n$, and thus
$q_\nu\=1-\rho(\bpnu)>t$.
Similarly, if $q<t<1$, then, for large $n$, $f_{X_\nu}(t)<t$  and thus
$q_\nu<t$. It follows that $q_\nu\to q$ as \ntoo.
\end{proof}

\begin{remark}   \label{Rgw}
 The case $\P(X=1)=1$, \ie{} $X=1$ \as, really is an exception.
If we let $X_\nu\sim\Be(1-\nu\qi)$, we have $X_\nu\dto X=1$, but
$\rho(\bpnu)=0$ for every $\nu$ while $\rho(\bp)=1$.
\end{remark}

\section{The giant component}\label{Spf}
Our ultimate goal is to describe the large component(s) of $\gndx\vvp$ 
and $\gnd\vvp$,
where $\vvq$ is one of the vaccination strategies defined above.
The basic strategy
will be to relate the neighbourhoods of a vertex
to a branching process.
We do this for \gndx, which is technically easier to handle;
as explained in \refSS{SSsimple}, the results
then transfer to \gnd{} too, provided \refCN{C2} holds.
We first do the argument in detail in the simplest case, \viz{}
$\gndx$ without edge deletion (i.e.\ $p=1$) or vaccination and
prove our main results
concerning the existence, size and uniqueness of the giant component.
We use and adapt the method in Bollob\'as, Janson and Riordan \cite{SJ178}
(for a different random graph model).
This will provide a new proof of the results by Molloy and Reed
\cite{MR1,MR2} (under our slightly weaker condition). We will then describe
the modifications needed to make the results valid also when there is
edge deletion or vaccination.

We say that an event holds \emph{with high probability} (\whp), if
it holds with probability tending to 1 as $n\to\infty$. We shall
use 
$o_p$
in the standard way (see \eg{} Janson, {\L}uczak and
Ruci\'nski~\cite{JLR}); for example,
if $(X_n)$ is a sequence of random variables, then
$X_n=o_p(1)$ means that $X_n \pto0$. We shall often use the
basic fact that, if $a\in\bbR$,
then $X_n\pto a$ if and only if, for every $\varepsilon >0$, the
relations $X_n>a-\eps$ and $X_n<a+\eps$ hold \whp{}. All unspecified
limits are taken as \ntoo, while $p$ and the vaccination
parameters $v$ or $c$ are kept fixed.

We denote the orders of the components of a graph $G$ by
$C_1(G)\ge C_2(G) \ge\dots$, with $C_j(G)=0$ if $G$ has fewer than
$j$ components. We let $\nk(G)$ denote the total number  of
vertices in components of order $k$, and write $\ngek(G)$ for
$\sum_{j\ge k}\nx{j}(G)$, the number of vertices in components of
order at least $k$.
Similarly,
we let $\ndk(G)$ and
$\ndgek(G)$ denote the number of such vertices that have degree $d$.

\begin{remark}
  Our results are typically of the form
$C_1(G_n)=\tau n+o_p(n)$ and $C_2(G_n)=o_p(n)$ for some number
  $\tau\ge0$
(or, equivalently,
$C_1(G_n)/n\pto\tau$ and $C_2(G_n)/n\pto0$).
Hence, if $\tau>0$, then there is exactly one ``giant''
component, and all other components are much smaller.
In our epidemic setting, this means that if $\tau=0$, then every
epidemic will be ``small'', \ie{} $o(n)$, while if $\tau>0$, then the
epidemic is large with probability $\tau$
(allowing the case that the initially infected person is vaccinated and thus
never becomes ill), and in that case, a
fraction $\tau$ of the population will be infected.
($\tau$ thus has a double role.)
\end{remark}

\subsection{$\gndx$, with $p=1$ and no vaccination}
As said above, we will use a branching process approximation. The
particles in the branching process correspond to free (not yet paired)
half-edges. Note that there are $jn_j$ half-edges belonging to
vertices of degree $j$. Hence, a random half-edge shares a vertex with
$j-1$ other half-edges with probability $jn_j/\sum_k kn_k$.
By \refCN{C1}, $jn_j/\sum_k kn_k\to jp_j/\mu$, and recall the
definition of $\tilde p_j=jp_j/\mu$ defined in (\ref{ptilde}).
Let $\bp$ be the Galton--Watson branching
process starting with one particle and with the offsping distribution
$(\tilde p_{j+1})_{j=0}^\infty$.
(This is the distribution $(p_j)_j$ size-biased and shifted one step.)
In other words, the offspring distribution is $\tD-1$, with $\tD$ as
in \refS{Sresults}. 

We let $\rho=\rho(\bp)$ denote the survival probability of $\bp$,
and define
\begin{equation}
  \tau\=\sumdoo p_d \bigpar{1-(1-\rho)^d};
\end{equation}
this is the survival probability for the branching process $\bp$
started with a random number of particles having the distribution
$(p_d)_{d=0}^\infty$.

Consider a vertex $x$ of degree $d$ in $\gndx$. We explore the component
containing $x$ by a breadth-first search.
We concentrate on the half-edges, so we begin by taking the $d$
half-edges at $x$, and label them as \emph{active}. We then process the active
half-edges one by one as follows. We take an active half-edge, relabel
it as \emph{used}, and find the half-edge that it connects to and the
corresponding vertex; this partner is chosen uniformly among
all half-edges that are not yet used.
We then label the partner as used and all other half-edges at the same
vertex as active, provided
that they are not already used (which would mean that we have found a
cycle or a multiple edge).
The active half-edges will behave essentially as a Galton--Watson
process (where we reveal the children of the particles one by one),
but the probability distribution of the children will vary slightly; it will
depend on the numbers of vertices of different degrees that we already
have found. Nevertheless, it is obvious that at each step in the beginning, the
probability of $j-1$ new half-edges is close to
$j n_{j}/ \sum_k k n_k \approx \tilde p_{j}$.

To be more precise,
first, let $k$ be a fixed number, and consider the event that $x$
belongs to a component with at least $k$ vertices. This is almost the
same as the probability that we will find at least $k-1$ active half-edges in
the process just described. (This is not exact, because if we stop
when we have found $k-1$ half-edges, some of these may connect back to
vertices already found; the probability of this tends to 0, however, as \ntoo.)
The complementary event, that the process finds less than $k-1$ active
half-edges, consists of a finite number of cases, where each case
describes the sequence of new active half-edges found at each step.
It is obvious that the probability of each of these cases converges,
as \ntoo, to the corresponding probability in $\bpd$, and thus
we find, for a vertex $x$ of degree $d$, with $\cC(x)$ denoting the
corresponding component of $\gndx$,
\begin{equation}\label{a1}
\P(|\cC(x)|\ge k)=\P(|\bpd|\ge k-1)+o(1).
\end{equation}
Recall that $\ngek$ is the number of vertices of degree $d$ belonging
to a component of size at least $k$. The expectation $\E\ndgek$ equals
$n_d$ times the probability that a given
vertex $x$ of degree $d$ satisfies $|\cC(x)|\ge k$, and thus, by
\eqref{a1} and \refCN{C1}(i), for every fixed $d\ge0$ and $k\ge1$,
\begin{equation}\label{a4}
\E(\ndgek /n) \to p_d\P(|\bpd|\ge k-1).
\end{equation}

We next want to let $k\to\infty$ here. We thus, for the remainder of
this section, assume that $\go(n)$ is a function such that
$\go(n)\to\infty$ but $\go(n)/n\to0$ as \ntoo. We regard components as
\emph{big} if they contain at least $\go(n)$ vertices, and
\emph{small} otherwise. (The flexibility in the choice of $\go(n)$ is
useful, but we will see that it does not matter much; the asymptotics
we find do not depend on $\go$.)

\begin{lemma}\label{L1}
If $\go(n)\to\infty$ and $\go(n)/n\to0$, then,
\begin{equation}\label{a2}
\E(\ngegon /n) \to \tau
\end{equation}
and,
for every fixed $d\ge0$,
\begin{equation}\label{a3}
\E(\ndgegon /n) \to p_d\P(|\bpd|=\infty)
=
p_d\bigpar{1-(1-\rho)^d}
.
\end{equation}
\end{lemma}

\begin{proof}
  We begin with an upper bound in \eqref{a3}.
For any fixed $k$, we have $\go(n)>k$ for large $n$, and thus
$\ndgegon\le\ndgek$. Consequently, \eqref{a4} yields
\begin{equation}\label{a5b}
\limsup_\ntoo\E(\ndgegon /n )
\le
\limsup_\ntoo\E(\ndgek /n )
= p_d\P(|\bpd|\ge k-1)
.
\end{equation}
As \ktoo, the right hand side converges to $p_d\P(|\bpd|=\infty)$,
and we find
\begin{equation}\label{a5}
\limsup_\ntoo\E(\ndgegon /n )
\le
 p_d\P(|\bpd|= \infty)
.
\end{equation}

For a lower bound, let $\nu\ge1$ be fixed and let
$X_\nu$ be a random variable
taking values in \set{0,1,\dots,\nu} with
$\P(X_\nu=j)=(1-\nu\qi)\tp_{j+1}$ for
$1\le j\le\nu$ (and a suitable value for $\P(X_\nu=0)$ so that the sum
becomes 1).
Consider the breadth-first exploration process described above.
As long as we have found less than $\go(n)$ vertices,
the number of new active half-edges at each step
stochastically dominates $X_\nu$, provided $n$ is large enough,
since the remaining number of vertices of degree $j+1$ is
$n_{j+1}-o(n)=p_{j+1} n-o(n)\ge(1-\nu\qi)p_{j+1} n$ for $n$ large.
(If $p_{j+1}=0$, the result is trivial.)
Consequently, letting $\bpnud$ be the Galton--Watson process with
$d$ initial particles and the
number of children distributed as $X_\nu$,
if  $n$ is large enough,
we can couple the exploration process and $\bpnud$ such that as long
as we have found less than $\go(n)$ vertices, the number of active
half-edges is at least the number of active particles in $\bpnud$
(\ie, the particles whose children have not yet been revealed.)
In particular, if the exploration process stops before $\go(n)$
vertices are found, then $\bpnud$ stops, and thus
the probability that a vertex $x$ of degree $d$ satisfies $|\cC(x)|< \go(n)$
is at most $\P(|\bpnud|<\infty)$.
Consequently, for large $n$,
\begin{equation}\label{a6}
\E\ndgegon \ge n_d \P(|\bpnud|=\infty)
\end{equation}
and thus
\begin{equation}
\liminf_\ntoo \E(\ndgegon /n )
\ge
p_d\P(|\bpnud|=\infty).
\end{equation}
Now let $\nu\to\infty$. Then $X_\nu\dto X$, where $X$ has the
distribution $\P(X=j)=\tilde p_{j+1}$, and thus, by \refL{Lgw},
$\P(|\bpnud|=\infty)\to\P(|\bpd|=\infty)$.
Consequently,
\begin{equation*}
\liminf_\ntoo \E(\ndgegon /n) \ge
p_d\P(|\bpd|=\infty),
\end{equation*}
which together with \eqref{a5} and \eqref{erika} yields \eqref{a3}.

Finally, noting that $\ndgegon \le n_d$, it follows easily from
the uniform summability in \eqref{unif1}
that we can sum \eqref{a3} over $d$ and take the limit outside the sum,
\ie{}
\begin{equation*}
\E(\ngegon /n)
=
\sum_d \E(\ndgegon /n)
\to
\sum_d p_d\bigpar{1-(1-\rho)^d}
=\tau.
\end{equation*}
\qedkorr
\end{proof}

Note that the limits do not depend on the choice of $\go(n)$. Hence,
it follows that the expected number of vertices belonging to
components of size between, say, $\log n$ and $n^{0.99}$ is $o(n)$.

We next show that we have convergence not only of the expectations but
also of the random variables in \eqref{a2} and \eqref{a3}, \ie{} that
these random variables are concentrated close to their expectations.

\begin{lemma}\label{L2}
If $\go(n)\to\infty$ and $\go(n)/n\to0$, then,
\begin{equation}\label{b2}
\ngegon /n \pto \tau
\end{equation}
and,
for every fixed $d\ge0$,
\begin{equation}\label{b3}
\ndgegon /n \pto
p_d\bigpar{1-(1-\rho)^d}
.
\end{equation}
\end{lemma}

\begin{proof}
Start with two distinct vertices $x$ and $y$ of the same degree $d$
and explore their components as above.
We can repeat the arguments above, and find
\begin{equation*}
\P\bigpar{|\cC(x)|< k,\,|\cC(y)|<k}=\P(|\bpd|< k-1)^2+o(1)
\end{equation*}
and thus, using \eqref{a1},
\begin{equation*}
\P\bigpar{|\cC(x)|\ge k,\,|\cC(y)|\ge k}=\P(|\bpd|\ge k-1)^2+o(1).
\end{equation*}
Multiplying with the number $n_d(n_d-1)$ of pairs $(x,y)$ of the same
degree $d$, and noting that the number of such pairs where both $x$
and $y$ belong to components of size $\ge k$ (the same or not) is
$\ndgek(\ndgek-1)$,
we find
\begin{equation*}
\E(\ndgek^2 /n^2)
=
\E(\ndgek(\ndgek-1)/n^2) +O(1/n)
\to p_d^2\P(|\bpd|\ge k-1)^2.
\end{equation*}
Hence,
$\limsup_\ntoo
\E(\ndgegon^2 /n^2)
\le p_d^2\P(|\bpd|\ge k-1)^2$ for every $k$, and thus
\begin{equation*}
\limsup_\ntoo \E(\ndgegon^2 /n^2)
\le p_d^2\P(|\bpd|=\infty)^2.
\end{equation*}
Since, by the \CS{} inequality and \eqref{a3}, further
\begin{equation*}
\E(\ndgegon^2 /n^2)
\ge (\E(\ndgegon /n))^2
\to p_d^2\P(|\bpd|=\infty)^2,
\end{equation*}
it follows that
\begin{equation*}
\E(\ndgegon^2/n^2) \to p_d^2\P(|\bpnud|=\infty)^2.
\end{equation*}
This and \eqref{a3} show that
\begin{equation*}
\Var\bigpar{\ndgegon/n} \to 0,
\end{equation*}
and thus
\begin{equation*}
\bigpar{\ndgegon-\E(\ndgegon)}/n \pto0,
\end{equation*}
which by \eqref {a3} implies \eqref{b3}.

Finally, again we can sum over $d$ because of \eqref{unif1}; this
yields \eqref{b2}.
\end{proof}

\begin{theorem}\label{T1}
Assume that Condition \ref{C1} holds.
Then
\begin{align*}
C_1(\gndx)&= \tau n + o_p(n),
\\
C_2(\gndx)&= o_p(n).
\end{align*}
\end{theorem}

\begin{proof}
We have already shown that roughly $\tau n$ vertices lie in big
components. It remains to show that most of them belong to the same
component. We write $G_n=\gndx$.

First, if $C_1(G_n)\ge \go(n)$, then
$\ngegon(G_n)\ge C_1(G_n)$. Thus, for every $\eps>0$ and $n$ so large
that $\go(n)<\eps n$, we have by \refL{L2}
\begin{equation}\label{jesper}
  \P\bigpar{C_1(G_n)>\tau n+\eps n}
\le   \P\bigpar{\ngegon(G_n)>\tau n+\eps n}
\to0.
\end{equation}

This completes the proof if $\tau=0$.

In the sequel we assume $\tau>0$ and show a corresponding estimate
from below.
First, if $p_d=0$ for every $d\ge2$, then $\tp_{j+1}=0$ for all $j\ge1$,
so $\bp$ dies immediately and $\rho=0$ and $\tau=0$.
Hence $p_d>0$ for some $d\ge2$. We fix such a $d$ for the
remainder of the proof, and fix $\gd$ with $0<\gd<1/2$.
Further, take (rather arbitrarily) $\go(n)=n^{0.9}$.

We assume in the sequel that $n$ is so large that
$n_d>\nigd$. We then
split the $\nigd$ first of the vertices of degree $d$ in $G_n$
into $d$ vertices of degree 1 each; we colour these $d\nigd$ new
vertices \emph{red}. (To be precise, we should round $\nigd$ to an integer.)
We denote the resulting graph by $G_n'$; note that $G_n'$ is a random
multigraph $G^*(n',(d_i'))$ where $n_j'$, the number of vertices of degree $j$,
is given by $n_d'=n_d-\nigd$, $n_1'=n_1+d\nigd$, and $n_j'=n_j$ for
$j\neq 1,d$. Note that the total number of vertices in $G_n'$ is
$n'\=n+(d-1)\nigd=n+o(n)$, and that $(d_i')$ satisfies \refCN{C1} with
the same $(p_j)$ (except that $n$ is replaced by $n'$, which only
makes a notational difference). Consequently, our results above apply
to $G_n'$ too.

By symmetry, we may assume that the $d\nigd$ red vertices in $G_n'$
are chosen at random among all vertices of degree 1, and that
$G_n$ is obtained by partitioning the red vertices at random into
groups with $d$ vertices and then coalescing each group into one vertex.

During the exploration of the component $\cC'(x)$ in $G'_n$ containing a vertex
$x$, in each step, the active half-edge is paired with the single
half-edge leading to a
red vertex with probability at least $c_1 \ngd$, for some $c_1>0$,
unless at least $\nigd$ red vertices already have been found.
Consequently, if the component $\cC'(x)$ has at least $\go(n)$ vertices, the
number of red vertices stochastically dominates
$\min\bigpar{\nigd,\Bi(\go(n)-1, c_1\ngd)}$.
A Chernoff bound, see \eg{} \cite[Corollary 2.3]{JLR},
shows that the probability that $\cC'(x)$ has at least $\go(n)$
vertices but less than $c_2\ngd \go(n)=c_2n^{0.9-\gd}$ red vertices
is at most $\exp(-c_3 n^{0.9-\gd})=o(n^{-1})$, for $c_2=c_1/2$ and
some $c_3>0$. Summing over all $x$,
we see that \whp, \emph{every} big component of $G'_n$
contains at least $c_2n^{0.9-\gd}$ red vertices.

Assume that this holds, and consider two big components $K_1$ and
$K_2$ in $G_n'$.
We can construct the random partition of the red vertices by taking
first the red vertices in $K_1$ one by one, unless already used, and
randomly selecting $d-1$ partners. We thus do this at least
$m:=c_2n^{0.9-\gd}/d$ times, and each time the probability of not
including a red vertex in $K_2$ is at most
$1-c_2n^{0.9-\gd}/(d\nigd)=1-c_4n^{-0.1}$, with $c_4=c_2/d$.
Consequently, the probability of not joining $K_1$ and $K_2$ in the
coalescing phase is at most
\begin{equation*}
 \exp(-mc_2 n^{-0.1})=\exp(-c_4^2 n^{0.8-\gd})
=o(n^{-2}).
\end{equation*}
Since there are at most $(n')^2=O(n^2)$ such pairs $K_1$ and $K_2$, we
see that \whp{} all big components in $G_n'$ are connected in $G_n$.
Hence, if $B'$ is the union of all big components in $G_n'$, and $B$
is the corresponding set of vertices in $G_n$, we see that \whp{} $B$
is connected in $G_n$, and, using \refL{L2} for $G_n'$,
\begin{equation}\label{emma}
  C_1(G_n) \ge |B|
\ge |B'|-(d-1)\nigd
=\tau n' +o_p(n)
=\tau n + o_p(n).
\end{equation}
Combining \eqref{emma} and \eqref{jesper} we obtain $C_1(G_n)=\tau n+o_p(n)$.

Finally, we observe that
if $C_2(G_n)\ge \go(n)$, then
$\ngegon(G_n)\ge C_1(G_n)+C_2(G_n)$, and thus,
by \eqref{b2} and \eqref{emma},
\begin{equation*}
  C_2(G_n) \le \max\bigpar{\go(n),\ngegon(G_n)-C_1(G_n)} =o_p(n).
\end{equation*}
\qedkorr
\end{proof}

\subsection{The simple random graph \gnd}\label{SSsimple}
We transfer the results to the simple random graph \gnd{} by the
following result proved in \cite{SJ195}; see also \eg{}
Bollob\'as \cite{bollobas} and McKay \cite{McKay} for
earlier versions.
\begin{lemma}
  \label{Lsimple}
If Conditions \refand{C1}{C2} hold, then
\begin{equation*}
\liminf_{\ntoo} \P\bigpar{\gndx\text{ is a simple graph}}>0.
\end{equation*}
\end{lemma}

All results for \gndx{} that can be stated in terms of convergence in
probability, as our results in this section, thus hold also if we
condition on the graph being simple.
In other words, the results proved for \gndx{} hold for \gnd{} too.
Thus, \refT{T1} has the following version for \gnd.
\begin{theorem}\label{T1s}
Assume that Conditions \refand{C1}{C2} hold.
Then
\begin{align*}
C_1(\gnd)&= \tau n + o_p(n),
\\
C_2(\gnd)&= o_p(n).
\end{align*}
\end{theorem}

\subsection{Uniform vaccination}
We now extend \refT{T1} to the graph $\gndx\uvp$ where $0\le v<1$ and
$0<p\le1$, see \refS{Smodels}.
Recall that we obtain this graph from $\gndx$ by randomly and
independently deleting edges with probability $1-p$ (non-transmission)
and vertices with probability $v$ (vaccination).
The branching process approximation arguments above still work, with
the difference that each new individual found is kept with probability
$p(1-v)$, and otherwise discarded.
Hence the offspring distribution is changed from $\tD-1$ to
$\tX_v\sim\MixBi(\tD-1,p(1-v))$, and the branching process
corresponding to an unvaccinated person with $d$ friends starts with
$\Bi(d,p(1-v))$
individuals. Let now $\bpd$ denote the branching process with this
offspring distribution, starting with $d$ individuals.
The probability generating function of $\tX_v$ is, as shown in
Subsections \refand{SSomodel}{SSumodel}, given by
\begin{equation*}
\E t^{\tX_v}= \frac{f_D'\bigpar{1-p(1-v)(1-t)}}{f_D'(1)}.
\end{equation*}
Hence, the extinction probability of $\bp^1$ is $\tpi\uvp$ given by
\eqref{tildepi_vU}. If we start the branching process with
$D'\sim\Bi(d,p(1-v))$ individuals, the extinction probability is thus,
writing $\xp=p(1-v)$,
\begin{equation*}
\pidu:=
  \sum_k \binom dk \xp^k(1-\xp)^{d-k}(\tpi\uvp)^k
=\bigpar{1-\xp+\xp\tpi\uvp}^d.
\end{equation*}
The arguments in the proofs of Lemmas \refand{L1}{L2} show, recalling
that each vertex has probability $1-v$ of being unvaccinated, that
\eqref{b3} holds in the form
\begin{equation*}\label{b3v}
\ndgegon /n \pto
p_d(1-v)\bigpar{1-\pidu},
\end{equation*}
for every fixed $d\ge0$, assuming
$\go(n)\to\infty$ and $\go(n)/n\to0$. Hence,
\begin{equation*}\label{b2v}
\ngegon /(n(1-v)) \pto
\sum_d p_d(1-\pidu)
=1-\sum_d p_d\pidu
=1-f_D(1-\xp+\xp\tpi\uvp).
\end{equation*}
This limit equals $\tau\uvp:=1-\pi\uvp$ with $\pi\uvp$ given
by \eqref{pi_vU}.

To extend \refT{T1}, it remains to show that there is only one very
large component. More precisely, we show again that, with
$\go(n)=n^{0.9}$, there is \whp{} only one big component. We argue as
for \refT{T1}, splitting some vertices of degree $d$ in $G_n=\gndx$
into $d$ red vertices of degree 1, calling the resulting graph
$G_n'$.

We vaccinate the vertices in $G_n'$ with probability $v$
each, independently; we then recombine the red vertices to vertices of
degree $d$ in $G_n$ and consider each such vertex as vaccinated if at
least one of its red parts in $G_n'$ is. This means that some vertices
in $G_n$ are vaccinated with probability larger than $v$, but this
does not hurt since the aim of the argument is to provide a lower
bound for $C_1$, the size of the largest component, and any extra
vaccinations can only decrease $C_1$.

By a Chernoff bound, there are \whp{} at least $(1-v)n^{1-\gd}$
unvaccinated red vertices, and it follows as before that \whp{} every
big component of $(G_n')\uvp$ contains at least $c_2 n^{0.9-\gd}$ red
vertices (although the value of $c_2$ may change). Given two big
components $K_1$ and $K_2$ it follows similarly as before that with
probability $1-o(n^{-2})$ there exists a vertex in $G_n$ that is split
into $d$ red vertices, of which at least one is in $K_1$, at least one
in $K_2$, and all are unvaccinated. The proof is completed as before.

Consequently, using also \refL{Lsimple}, we have the following
theorem. Theorems \refand{TU}{T0} (the special case $v=0$) are
immediate consequences.

\begin{theorem}
  \label{T1U}
Assume \refCN{C1}, and let $0<p\le1$, $0\le v<1$.
Then,
\begin{align*}
C_1\bigpar{\gndx\uvp}&= \tau\uvp n(1-v) + o_p(n),
\\
C_2\bigpar{\gndx\uvp}&=  o_p(n),
\end{align*}
where $\tau\uvp=1-\pi\uvp$ with $\pi\uvp$ given by \eqref{pi_vU}.
If also \refCN{C2}  holds, then the same results hold for $\gnd\uvp$ too.
\end{theorem}

\subsection{Acquaintance vaccination}

As explaind in \refSS{SSamodel}, in order to obtain (asymptotically)
a Galton--Watson branching process, with the right independence
properties, we consider directed edges, or equivalently
half-edges, that are \emph{open}, \ie{} transmission may take place
but the edge is not used for vaccination.
Moreover, we consider only open edges originating at an
unvaccinated person.

Let $x$ be a given vertex with degree $d$ in \gndx, and let us explore
the component of $x$ in \gndxa, conditioned on $x$ being unvaccinated
(otherwise $x$ does not belong to \gndxa).
In order to be kept in \gndxa, an edge has to be open, but not all
edges are kept since some may lead to vertices that are vaccinated,
see Figure \ref{acq-br}c). Nevertheless, we consider all open edges
found during the exploration. We declare the open edges starting at
$x$ to be \emph{active}. We then investigate the active edges. If an
active edge leads to a person that is unvaccinated, we declare the
open edges going from that person, except the one going back to where
we just came from, to be new active edges.
We continue until no more active edges are found; we then have found
the component containg $x$ (plus some extra open edges leading to
vaccinated persons).

We investigate this process probabilistically, revealing the structure
of $\gndx$ by combining half-edges at random as we proceed the
exploration.
We consider asymptotics as \ntoo, and some of the statements below are only
approximatively correct for finite $n$.

Note first that each of the $d$ edges leading from $x$ is open with
probability $pe^{-c/d}$, independently of each other, so we start with
$\Bi(d,pe^{-c/d})$ open edges.

The vertex $x$ has $d$ friends; in \gndx{} they are chosen by randomly
choosing $d$ half-edges and their degrees have the size-biased
distribution $(\tp_j)$, independently of each other.
Conditioning on $x$ being unvaccinated means that we condition on none
of the $d$ edges being used for vaccination in the opposite
direction. Since the probability that a friend with degree $j$ does
not name $x$ is $e^{-c/j}$, this preserves the independence of the
degrees of the friends, but shifts their distribution to, as asserted
in \eqref{pfriend}, $(\tp_je^{-c/j}/\ga)_j$, where
$\ga=\ga(c)=\sum_j\tp_je^{-c/j}$ as in \eqref{alpha} is the probability of not
being named by a random friend.

Now suppose that an open edge goes from $x$ to a friend $y$ of degree
$k$. In order for this to define an edge in \gndxa, $y$ must not be
vaccinated through another of its friends; this has the probability $\ga^{k-1}$.
In this case, $y$ has $k-1$ further edges, and each of them is open
with probability $pe^{-c/k}$. It follows that the number of new open
edges at $y$ has a distribution that is the mixture
$(1-\ga^{k-1})\gd_0+\ga^{k-1}\Bi(k-1,pe^{-c/k})$.
Using the distribution \eqref{pfriend} for the degree of $y$, we
finally see that the distribution of the number $Y$ of new active
edges found when exploring a single active edge is given by \eqref{aoff}.

Hence, observing obvious independence properties, the process of
active edges is (asymptotically) a Galton-Watson branching process
with offspring distribution $Y$, starting with $\Bi(d,pe^{-c/d})$
active edges. Denote his branching process by $\bpxd$.
Let, as in \refSS{SSamodel}, $\tpi\acp$ by the probability that a
branching proess with this offspring distribution $Y$ and starting with a
single individual dies out. Then, the extinction probability of
$\bpxd$ is
\begin{equation*}
  \begin{split}
  \pida&:=\P(|\bpxd|<\infty)
= \sum_{j=0}^d\binom dj(pe^{-c/d})^j(1-pe^{-c/d})^{d-j}(\tpi\acp)^j
\\&
=
(1-pe^{-c/d}+pe^{-c/d}\tpi\acp)^d.
  \end{split}
\end{equation*}

A minor complication is that the branching process approximation
counts open edges and, as remarked above, not all open edges lead to
vertices in \gndxa.
Thus \eqref{a1} does not extend directly. However, we still have the inequality
\begin{equation*}
\P(|\cC(x)|\ge k)\le\P(|\bpxd|\ge k-1)+o(1).
\end{equation*}
Furthermore, a vertex of degree $d$ in $\gndx$ is unvaccinated with
probabilty $\ga^d$, and thus
\begin{equation*}
\E(\ndgek)\le n_d \ga^p\bigpar{\P(|\bpxd|\ge k-1)+o(1)},
\end{equation*}
which arguing as in \eqref{a5b} and \eqref{a5} leads to
\begin{equation}\label{f1}
\limsup_\ntoo\E(\ndgegon/n)
\le p_d \ga^p\P(|\bpxd|=\infty)
= p_d \ga^p(1-\pida).
\end{equation}

For a lower bound, we note that an open edge creates new open edges in
the exploration process only if it leads to an unvaccinated person.
Hence, if $f(\bpxd)$ denotes the number of individuals in the branching
process $\bpxd$ with at least one child, we have, for every $k\ge1$,
\begin{equation*}
\P(|\cC(x)|\ge k)\ge\P(f(\bpxd)\ge k-1)+o(1).
\end{equation*}
In order to replace the fixed $k$ by $\go(n)$, we do as in the proof
of \refL{L1} and define a Galton--Watson process $\bpxdnu$, now starting
with
$\Bi(d,pe^{-c/d}(1-\nu\qi))$  individuals
and with an offspring distribution $Y_\nu$ on \set{0,\dots,\nu} with
$\P(Y_\nu=j)=(1-\nu\qi)\P(Y=j)$ for $j=1,\dots,\nu$.

For each $\nu$ and each fixed $A<\infty$, we can for large $n$ couple
the exploration process and $\bpxdnu$ as in the proof of \refL{L1} as
long as we have found at most $A\go(n)$ open edges.
Hence, if $|\cC(x)|<\go(n)$, then either $f(\bpxdnu)<\go(n)$ or the
process $\bpxdnu$ reaches more than $A\go(n)$ individuals while less
than $\go(n)$ of them, plus the root, have had children. The
probability of the latter event is at most, since the root has at most
$d$ children,
$$
\P\Bigpar{1+d+\sum_{i=1}^{\go(n)} Y^*_{\nu,i} > A\go(n)},
$$
where $Y^*_{\nu,i}$ are independent random variables with the
distribution $\cL(Y\mid Y>0)$, and thus this probability tends to 0 by
the law of large numbers provided we have chosen $A>\E(Y\mid Y>0)$.

Consequently,
\begin{equation*}
  \P\bigpar{|\cC(x)|<\go(n)}
\le
 \P\bigpar{f(\bpxdnu)<\go(n)} + o(1)
\le
 \P\bigpar{f(\bpxdnu)<\infty} + o(1).
\end{equation*}
Using again that a person with degree $d$ is unvaccinated with
probability $\ga^d$, it follows that
\begin{equation*}
\E\ndgegon
\ge n_d \ga^p\bigpar{1-\P(|\bpxdnu|<\infty)+o(1)}
\end{equation*}
and thus
\begin{equation*}
\liminf_\ntoo\E(\ndgegon/n)
\ge p_d \ga^p\P(|\bpxdnu|=\infty).
\end{equation*}
We let $\nu\to\infty$ and obtain by \refL{Lgw}
\begin{equation*}
\liminf_\ntoo\E(\ndgegon/n)
\ge p_d \ga^p\P(|\bpxd|=\infty)
= p_d \ga^p(1-\pida),
\end{equation*}
which together with \eqref{f1} yields
\begin{equation*}
\E(\ndgegon/n)
\to p_d \ga^p(1-\pida).
\end{equation*}
Arguing as in the proof of \refL{L2}, we find also
\begin{equation*}
\ndgegon/n
\pto p_d \ga^p(1-\pida)
\end{equation*}
and, recalling \eqref{pi_vA} and \eqref{v_A},
\begin{equation*}
\ngegon/n
\pto \sum_d p_d \ga^p(1-\pida)
= \sum_d p_d \ga^p(1-\pi\acp)
= (1-v(c))\tau\acp,
\end{equation*}
with $\tau\acp=1-\pi\acp$. In particular,
\begin{equation*}
C_1(\gndxa)
\le\go(n)+\ngegon
\le (1-v(c))\tau\acp n+o_p(n).
\end{equation*}

Finally, we argue again as in the proof of \refT{T1} to show that most
vertices in large components belong to a single component. We split
some of the vertices in $G_n=\gndx$ as above and perform acquaintance
vaccination on the resulting graph $G'_n$. This corresponds to
acquaintance vaccination on $G_n$, except that the vertices that are
split now are asked to name a friend $\Po(dc)$ times instead of
$\Po(c)$. We perform thus some extra vaccinations, but this can only
decrease $C_1$ and we obtain as in \eqref{emma} the lower bound
\begin{equation*}
  C_1(\gndxa)
\ge (1-v(c))\tau\acp n + o_p(n).
\end{equation*}

Summing up, and using \refL{Lsimple}, we have shown he following
theorem.  \refT{TA} is an immediate consequence.

\begin{theorem}
  \label{T1A}
Assume \refCN{C1}, and let $0<p\le1$, $0\le c<\infty$.
Then, 
\begin{align*}
C_1\bigpar{\gndx\acp}&= \tau\acp n(1-v(c)) + o_p(n),
\\
C_2\bigpar{\gndx\acp}&=  o_p(n),
\end{align*}
where $\tau\acp=1-\pi\acp$ with $\pi\acp$ given by \eqref{pi_vA}.
If also \refCN{C2}  holds, then the same results hold for $\gnd\acp$ too.
\end{theorem}

\subsection{Edgewise vaccination}

We argue as for acquaintance vaccination with the modifications
(simplifications) explained in \refSS{SSemodel}. There are no new
complications, and we obtain the following.
\refT{TE} is an immediate consequence.

\begin{theorem}
  \label{T1E}
Assume \refCN{C1}, and let $0<p\le1$, $0<\ga\le1$.
Then, for $j=1,2$,
\begin{align*}
C_1\bigpar{\gndx\ejqp}&= \tau\ejqp n(1-v(\ga)) + o_p(n),
\\
C_2\bigpar{\gndx\ejqp}&=  o_p(n),
\end{align*}
where $\tau\ejqp=1-\pi\ejqp$ with $\pi\ejqp$ given by \eqref{pi_vE}.
If also \refCN{C2}  holds, then the same results hold for $\gnd\ejqp$ too.
\end{theorem}

\section*{Acknowledgements}
We thank Mathias Lindholm for help in producing the figures. T.B.\
gratefully acknowledges financial support from the Swedish Research Council.

\newcommand\AAP{\emph{Adv. Appl. Probab.} }
\newcommand\JAP{\emph{J. Appl. Probab.} }
\newcommand\AMS{Amer. Math. Soc.}
\newcommand\JAMS{\emph{J. \AMS} }
\newcommand\MAMS{\emph{Memoirs \AMS} }
\newcommand\PAMS{\emph{Proc. \AMS} }
\newcommand\TAMS{\emph{Trans. \AMS} }
\newcommand\AnnMS{\emph{Ann. Math. Statist.} }
\newcommand\AnnPr{\emph{Ann. Probab.} }
\newcommand\AnnSt{\emph{Ann. Statist.} }
\newcommand\AnnAP{\emph{Ann. Appl. Probab.} }
\newcommand\CPC{\emph{Combin. Probab. Comput.} }
\newcommand\JMAA{\emph{J. Math. Anal. Appl.} }
\newcommand\JASA{\emph{J. Amer. Statist. Assoc.} }
\newcommand\RSA{\emph{Random Struct. Alg.} }
\newcommand\ZW{\emph{Z. Wahrsch. Verw. Gebiete} }
\newcommand\LNCS[2]{Lecture Notes Comp. Sci. #1, Springer, #2}
\newcommand\Springer{Springer}
\newcommand\Wiley{Wiley}
\newcommand\CUP{Cambridge Univ. Press}
\newcommand\Poznan{Pozna\'n}
\newcommand\jour{\emph}
\newcommand\book{\emph}
\newcommand\inbook{\emph}
\newcommand\vol{\textbf}
\newcommand\toappear{\unskip, to appear}
\newcommand\webcitesvante{\webcite{http://www.math.uu.se/\~{}svante/papers}}
\newcommand\arxiv[1]{\webcite{http://arxiv.org/#1}}
\renewcommand\&{and}


\begin{thebibliography}{99}
\bibitem{SJ195}
S. Janson,
The probability that a random multigraph is simple.
\arxiv{math.CO/0609802}
\end{thebibliography}
\end{document}